\documentclass[12pt]{article}

\usepackage{amsmath,amssymb,enumerate,bbm}
\newtheorem{theorem}{Theorem}[section]

\newtheorem{remark}[theorem]{Remark}

\newtheorem{conjecture}[theorem]{Conjecture}
\numberwithin{equation}{section}

 \makeatletter
 
 \newcommand{\Rmnum}[1]{\expandafter\@slowromancap\romannumeral #1@}
 \makeatother

\begin{document}
\title{The Erd\"os-Falconer distance problem on the unit sphere in vector spaces over finite fields}\author{Le Anh Vinh\\
Mathematics Department\\
Harvard University\\
Cambridge, MA 02138, US\\
vinh@math.harvard.edu}\maketitle

\begin{abstract}
Hart, Iosevich, Koh and Rudnev (2007) show, using Fourier analysis method, that the finite Erd\"os-Falconer distance conjecture holds for subsets of the unit sphere in $\mathbbm{F}_q^d$. In this note, we give a graph theoretic proof of this result.  
\end{abstract}

\section{Introduction}

The Erd\"os Distance Problem is perhaps the best known problem in combinatorial
geometry. How many distinct distances can occur among $n$ points in the plane?
Although this problem has received considerable attention, we are still far
from the solution. The Falconer distance conjecture says that if $E \subset \mathbbm{R}^d$, $d \geq 2$, has Hausdroff dimension greater than $\frac{d}{2}$, then the set of distances occur in $E$ has positive Lebesgue measure. See \cite{iosevich-rudnev-uriarte} for the connections between the Erd\"os and Falconer distance conjectures.

In the finite field setting, the distance problem turns out to have features of both the Erd\"os and Falconer distance problems in real spaces. Let $\mathbbm{F}_q$ denote the finite field with
$q$ elements where $q \gg 1$ is an odd prime power. For any $x, y \in
\mathbbm{F}_q^d$, the distance between $x, y$ is defined as $\|x - y\|= (x_1 -
y_1)^2 + \ldots + (x_d - y_d)^2$. Let $E \subset \mathbbm{F}_q^d$, $d
\geqslant 2$. Then the finite analog of the classical Erd\"os distance problem is
to determine the smallest possible cardinality of the set
\begin{equation}
  \Delta (E) =\{\|x - y\|: x, y \in E\},
\end{equation}
viewed as a subset of $\mathbbm{F}_q$. Bourgain, Katz and Tao (\cite{bourgain-katz-tao}), showed,
using intricate incidence geometry, that for every $\varepsilon > 0$, there
exists $\delta > 0$, such that if $E \in \mathbbm{F}_q^2$ and
$C_{\varepsilon}^1 q^{\varepsilon} \leqslant |E| \leqslant C_{\varepsilon}^2
q^{2 - \varepsilon}$, then $| \Delta (E) | \geqslant C_{\delta} |E|^{\frac{1}{2}
+ \delta}$ for some constants $C_{\varepsilon}^1, C_{\varepsilon}^2$ and
$C_{\delta}$. The relationship between $\varepsilon$ and $\delta$ in their
argument is difficult to determine. Going up to higher dimension using
arguments of Bourgain, Katz and Tao is quite subtle. Iosevich and Rudnev \cite{iosevich-rudnev}
establish the following result using Fourier analytic method.

\begin{theorem}(\cite{iosevich-rudnev}) \label{ir1}
  Let $E \subset \mathbbm{F}_q^d$ such that $|E| \gtrsim C q^{d / 2}$ for
  $C$ sufficiently large. Then
  \begin{equation}
    | \Delta (E) | \gtrsim \min \left\{ q, \frac{|E|}{q^{(d - 1) / 2}}
    \right\} .
  \end{equation}
\end{theorem}

In view of this reslut, Iosevich and Rudnev (\cite{iosevich-rudnev}) formulated the Erd\"os-Falconer conjecture as follows.

\begin{conjecture} \label{efc}
Let $E \subset \mathbbm{F}_q^d$ such that $|E| \geq C_{\epsilon} q^{\frac{d}{2} + \epsilon}$. Then there exists $c > 0$ such that $|\Delta(E)| \geq cq$.
\end{conjecture}

By modifying the proof of Theorem \ref{ir1} slightly, Iosevich and Rudnev (\cite{iosevich-rudnev}) obtain the following stronger conclusion.

\begin{theorem}(\cite{iosevich-rudnev})\label{ir2}
Let $E \subset \mathbbm{F}_q^d$ such that $|E| \geq Cq^{\frac{d+1}{2}}$ for sufficient large constant $C$. Then $\Delta(E) = \mathbbm{F}_q$. 
\end{theorem}

In \cite{hart-iosevich-koh-rudnev}, the authors show that Theorem \ref{ir2} is essentially sharp, which implies that Conjecture \ref{efc} is not true in general. They show however, that the exponent predicted by Conjecture \ref{efc} does hold for subsets of the sphere $S^{d-1} = \{x\in \mathbbm{F}_q^d : x_1^2 + \ldots x_d^2 = 1\}$.

\begin{theorem} (\cite{hart-iosevich-koh-rudnev}) \label{mt-dps}
Let $E \subset \mathbbm{F}_q^d$, $d\geq 3$, be a subset of the sphere \[S^{d-1} = \{x\in \mathbbm{F}_q^d: ||x|| = 1\}.\] Suppose that $|E| \geq Cq^{d/2}$ with a sufficiently large constant $C$. Then there exists $c > 0$ such that $|\Delta(E)| > cq$. 
\end{theorem}

In this note, we will give a graph theoretic proof of this result. The rest of this note is organized as follows. In Section 2, we construct our main tools to study the Erd\"os-Falconer distance problem over subsets of the sphere, the graphs associated to the projective spaces over finite fields. Our construction follows one of Bannai, Shimabukuro and Tanaka in \cite{bannai-shimabukuro-tanaka}. We then prove Theorem \ref{mt-dps} in Section 3. We also call the reader's attention that this note is a subsequent of an earlier paper \cite{vinh-ejc1}.

\begin{remark} If $d$ is even, Hart, Iosevich, Koh and Rudnev showed that all the ditances can be obtained under the same assumption and the size condition on $E$ cannot be relaxed. If $d$ is odd then we cannot in general get all the distances if $|E| \ll q^{\frac{d+1}{2}}$. Interested readers can see \cite{hart-iosevich-koh-rudnev} for a detailed discussion and related results. 
\end{remark}

\section{Finite non-Euclidean graphs}

In this section, we give a construction of graphs from the action of simple orthogonal group on the set of non-isotropic square-type of projective spaces over finite fields. Our construction follows one of Bannai, Shimabukuro and Tanaka in \cite{bannai-shimabukuro-tanaka}. Let $V = \mathbbm{F}_q^d$ be the $d$-dimensional vector space over the finite field $\mathbbm{F}_q$ ($q$ is an odd prime power). For each element $x$ of $V$, we denote the $1$-dimensional subspace containing $x$ by $[x]$. Let $\Omega$ be the set of all square type non-isotropic $1$-dimensional subspaces of $V$ with respect to the quadratic form $Q(x) = x_1^2 + \ldots + x_d^2$. The simple orthogonal group $O_d(\mathbbm{F}_q)$ acts transtively on $\Omega$, and yields a symmetric association scheme $\Psi(O_d(\mathbbm{F}_q),\Omega)$ of class $(q+1)/2$. We have two cases.

Case I. Suppose that $d = 2m+1$. The relations of $\Psi(O_{2m+1}(\mathbbm{F}_q),\Omega)$ are given by
\begin{eqnarray*}
 R_1 & = & \{([U],[V]) \in \Omega \times \Omega \mid  (U+V) \cdot (U+V) = 0\},\\
 R_i & = & \{([U],[V]) \in \Omega \times \Omega \mid (U+V) \cdot (U+V) = 2 + 2 \nu^{- (i
  - 1)} \} \, (2 \leqslant i \leqslant (q - 1) / 2)\\
 R_{(q+1)/2} & = & \{([U], [V]) \in \Omega \times \Omega \cdot (U+V) \cdot (U+V) = 2\},
\end{eqnarray*}
where $\nu$ is a generator of the field $\mathbbm{F}_q$ and we assume $U\cdot U = 1$ for all $[U] \in \Omega$ (see \cite{bannai-hao-song}).

Case II. Suppose that $d = 2m$. The relations of $\Psi(O_{2m}(\mathbbm{F}_q),\Omega)$ are given by
\begin{eqnarray*}
 R_i & = & \{([U],[V]) \in \Omega \times \Omega \mid (U+V) \cdot (U+V) = 2 + 2^{-1} \nu^{i} \} \, (1 \leqslant i \leqslant (q - 1) / 2)\\
 R_{(q+1)/2} & = & \{([U], [V]) \in \Omega \times \Omega \cdot (U+V) \cdot (U+V) = 2\},
\end{eqnarray*}
where $\nu$ is a generator of the field $\mathbbm{F}_q$ and we assume $U\cdot U = 1$ for all $[U] \in \Omega$ (see \cite{bannai-hao-song}). 

The graphs $(\Omega,R_i)$ are not Ramanujan in general, but fortunately, they are asymptotic Ramanujan for large $q$. The following theorem summaries the results from \cite{bannai-shimabukuro-tanaka} in a rough form.

\begin{theorem}(\cite{bannai-shimabukuro-tanaka})\label{mt-bst} The graphs $(\Omega,R_i)$ $(1 \leq i \leq (q+1)/2)$ are regular of order $q^{d-1}(1+o(1))/2$ and valency $Kq^{d-2}(1+o(1))$. Let $\lambda$ be any eigenvalue of the graph $(\Omega,R_i)$ with $\lambda \neq$ valency of the graph then 
\[|\lambda| \leq k(1+o(1))q^{(d-2)/2},\]
for some $k,K > 0$ (In fact, we can show that $k = 2$ and $K = 1$ or $1/2$).
\end{theorem}

\section{Graph theoretic proof of Theorem \ref{mt-dps}}

Let $E$ be a subset of the unit sphere $S^{d-1} = \{x\in \mathbbm{F}_q^d: ||x|| = 1\}$ with $|E| \geq Cq^{d/2}$. Let $E_1 = \{[x] : x \in E\} \subset \Omega$ (where $\Omega$ is the set of all square type non-isotropic $1$-dimensional subspaces of $V$ with respect to the quadratic form $Q(x) = x_1^2 + \ldots + x_d^2$). Since each line through origin in $\mathbbm{F}_q^d$ intersects the unit sphere $S^{d-1}$ at two points, $|E_1| \geq Cq^{d/2}/2$. Suppose that $([U], [V]) \in E_1 \times E_1$ is an edge of $(\Omega,R_i)$. Then
\[(U+V)\cdot(U+V) = 2+\alpha_i,\]
where $\alpha_i = 2\nu^{(-(i-1))}$ if $d$ is odd and $\alpha_i = 2^{-1}\nu^i$ if $d$ is even. Since $U\cdot U = V \cdot V = 1$, we have $(U-V)\cdot (U-V) = 2 - \alpha_i$. The distance between $U$ and $V$ (in $E$) is either $(U+V)\cdot (U+V)$ or $(U-V) \cdot (U+V)$, so
\[|\Delta(E) \cap \{2+\alpha_i, 2-\alpha_i\}| \geq 1.\]
Therefore, it is sufficient to show that $E_1 \times E_1$ contains edges of at least $cq$ graphs among $(\Omega,R_i)$, $1\leq i \leq (q+1)/2$.

To complete the proof, we need the following result from spectral graph theory. We call a graph $G$ $(n, d,\lambda)$-regular if $G$ is a $d$-regular graph on $n$ vertices with the absolute value of each of its eigenvalues but the
largest one is at most $\lambda$. It is well-known that if $\lambda \ll d$ then a $(n,d,\lambda)$-regular graph behaves similarly as a random graph $G_{n,d/n}$. Presicely, we have the following result (see Corollary 9.2.5 and Corollary 9.2.6 in \cite{alon-spencer}).

\begin{theorem} \label{expander} (\cite{alon-spencer})
  Let $G$ be a $(n, d, \lambda)$-regular graph.  For every set of vertices $B$ of $G$, we have
  \begin{equation}\label{f2}
    |e_G(B) - \frac{d}{2 n} |B|^2 | \leqslant \frac{1}{2} \lambda |B|,
  \end{equation}
  where $e_G(B)$ is number of edges in the induced subgraph of $G$ on $B$.
\end{theorem}  

From Theorem \ref{expander}, we have
\[ |e_{(\Omega,R_i)}(E_1) - \frac{Kq^{d-2}(1+o(1))}{q^{d-1}(1+o(1))/2}|E_1|^2| \leq \frac{1}{2}k(1+o(1))q^{(d-2)/2}|E_1|.
\]
Since $|E_1| \geq Cq^{d/2}/2$ for $C$ sufficiently large, the left hand term $\frac{1}{2}k(1+o(1))q^{(d-2)/2}|E_1|$ is neglected by $\frac{Kq^{d-2}(1+o(1))}{q^{d-1}(1+o(1))/2}|E_1|^2$. 
Thus, we have $e_{(\Omega,R_i)}(E_1) = O(|E_1|^2/q)$. Besides $E_1 \times E_1$ is edge-decomposed into $(q+1)/2$ graphs $(\Omega,R_i)$, $1\leq i \leq (q+1)/2$. This implies that $E_1 \times E_1$ contains edges of at lesat $cq$ graphs among $(\Omega,R_i)$, $1\leq i \leq (q+1)/2$ for some constant $c > 0$. The theorem follows.


\begin{thebibliography}{00}




\bibitem{alon-spencer} N. Alon and J. H. Spencer,
\textit{The probabilistic method}, 2nd ed., Willey-Interscience, 2000.

\bibitem{bannai-hao-song} E. Bannai, S. Hao and S.-Y. Song, Character tables of the association schemes of finite orthogonal groups acting on the nonisotropic points, \textit{Journal of Combinatorial Theory, Series A} \textbf{54} (1990), 164-170.

\bibitem{bannai-shimabukuro-tanaka} E. Bannai, O. Shimabukuro and H. Tanaka, Finite analogues of non-Euclidean spaces and Ramanujan graphs, \textit{European Journal of Combinatorics} \textbf{25} (2004), 243--259.

\bibitem{bourgain-katz-tao}
J. Bourgain, N. Katz, T. Tao, A sum-product estimate in finite fields, and applications, \textit{Geom. Funct. Anal.} \textbf{14} (2004), 27-57.

\bibitem{brass}
P. Brass, W. Moser and J. Pach, \textit{Research problems in discrete geometry}, Springer, 2005.


\bibitem{erdos46} P. Erd\"os, On sets of distances of $n$ points, \textit{Amer. Math. Monthly} \textbf{53} (1946) 248--250.

\bibitem{hart-iosevich-koh-rudnev}D. Hart, A. Iosevich, D. Koh and M. Rudnev, Averages over hyperplanes, sum-product theory in vector spaces over finite fields and the
Erd\"os-Falconer distance conjecture, preprint, 2007.

\bibitem{iosevich-rudnev}
A. Iosevich, M. Rudnev, Erd\"os distance problem in vector spaces over finite fields, \textit{Transactions of the American Mathematical Society} \textbf{359} (12) (2007), 6127-6142.

\bibitem{iosevich-rudnev-uriarte} A. Iosevich, M. Rudnev and I. Uriarte-Tuero. Theory of dimension for large discrete sets and applications, preprint, 2007.

\bibitem{szekely97} L. A. Sz\'ekely, Crossing numbers and hard Erd\"os problems in discrete geometry, \textit{Comb. Probab. Comput.} \textbf{6} (1997), 353--358.

\bibitem{vinh-ejc1} L. A. Vinh, Explicit Ramsey graphs and Erd\"os distance problem over finite Euclidean and non-Euclidean spaces, \textit{Electronic Journal of Combinatorics} \textbf{15} (2008), R5.

\end{thebibliography}
\end{document}